\documentclass{article}
\usepackage[authoryear,sectionbib,sort]{natbib}
\usepackage{amsmath}
\usepackage{amssymb}
\usepackage{float}
\usepackage{graphicx}
\usepackage{cite}
\usepackage[doublespacing]{setspace}
\usepackage{hyperref}
\title{Large Deviations of the Estimated Cumulative Hazard Rate}
\author{Niklas Hohmann\thanks{GeoZentrum Nordbayern, Fachgruppe Pal\"aoumwelt, University of Erlangen-Nuremberg, Loewenichstr. 28, 91054 Erlangen, Germany\; email:\href{mailto:niklas.hohmann@fau.de}{niklas.hohmann@fau.de}}}
\begin{document}
\maketitle
\begin{abstract}
Survivorship analysis allows to statistically analyze situations that can be modeled as waiting times to an event. These waiting times are characterized by  the cumulative hazard rate, which can be estimated by the Nelson-Aalen estimator or diverse confidence estimators based on asymptotic statistics. To better understand the small sample properties of these estimators, the speed of convergence of the estimate to the exact value is examined. This is done by deriving large deviation principles and their rate functions for the estimators and examining their properties. It is shown that these rate functions are asymmetric, leading to a tendency of the estimated cumulative hazard rate to overestimate the true cumulative hazard rate. This tendency is strongest in the cases of (1) small sample sizes and (2) low tail probabilities. Taking this tendency into account can improve risk assessments of rare events and of cases where only little data is available.
\end{abstract}


\section{Introduction}
\label{sec:Introduction}
Survival analysis is the standard framework to statistically analyze waiting times until an event occurs \citep{milicic2008}. Classical applications of survival analysis arise in epidemiology, where the event can be the recovery or the death of a patient, and the results of the statistical analysis can decide on the admission or nonadmission of a new medical procedure \citep{Sasieni2014,kantoff2010}.\\
A fundamental concept in characterizing the waiting time up to an event is the cumulative hazard rate (short CH)\citep{aalen2008}. It is estimated using the Nelson-Aalen estimator and can be complemented by different confidence intervals and bands to compensate for uncertainties of the results \citep{aalen2008,aalen1978,bie1986,nelson1969,nelson1972}.\\
However these confidence area estimators are based on the asymptotic properties of the Nelson-Aalen estimator\citep{bie1986}, making it impossible to derive analytical statements regarding their performance for small sample sizes. As a result, no reliable assessments of the uncertainties for small samples are available, which is problematic for studies where larger sample sizes cannot be achieved, be it because of ethical considerations, financial reasons, or rarity.\\
The aim of this paper is to derive analytical expressions regarding the small sample properties of the estimated CH to better assess its deviations from the true CH in these cases.\\
For this, the pointwise speed of convergence of the estimated CH to the true CH is examined using the theory of large deviations.

\section{Outline}
\label{sec:Outline}
The conventions $\ln(0):= \lim_{x \downarrow 0} \ln(x)= - \infty$ and $0 \ln(0):=0$ will be used throughout the paper. \\
Let $X, X_1, X_2, \dots$ be i.i.d. positive random variables, whose values model the waiting times for an event to occur. With the survival function $S(t):= P(X >t)$\citep[p. 9]{kleinbaum2010}, the cumulative hazard rate (short CH) can be written as\citep[p. 294]{kleinbaum2010}
\begin{equation}
H(t):=-\ln (S(t)).
\end{equation}
It can be estimated using the $n$-th empirical CH
\begin{equation}
H_n(t):= -\ln(S_n(t)),
\end{equation}
where $S_n(t)= \tfrac 1n \sum_{i=1}^n \mathbf 1_{(t,+ \infty)}(X_i)$ is the $n$-th empirical survival function. The aim of this paper is to examine the pointwise speed of convergence of $H_n$ towards $H$.\\
From the theory of large deviations, it is known that for integrable i.i.d. $(X_i)_{i \in \mathbb N}$, the relation
\begin{equation}
P\left( \tfrac 1n \sum_{i=1}^n X_i  \in C \right) \approx \exp(- n \inf I(C)) \label{eq:ldpfundamental}
\end{equation}
holds for a large class of sets $C$\citep{varadhan2016}. The function $I$ is called a rate function and determines the speed of convergence of the averaged random variables towards their expectation value.\\
In this paper, a rate function for the $H_n$ is derived (section \ref{sec:Establishing}). The properties of this rate function are then examined to draw conclusions about the convergence of $H_n$ to $H$ (section (\ref{sec:Properties}).\\
For this, fix any $\tilde t \in \mathbb{R}^+$, and define the tail probabilities as $p:=p(\tilde t)=P(X>\tilde t)$. Assume that $p \in (0,1)$, since the cases $p=0$ and $p=1$ are trivial. Define the i.i.d. random variables
\begin{equation}
Y_i= \begin{cases}1 & \text{ if } X_i > \tilde t \\ 0 & \text{ if } X_i \leq \tilde t \end{cases} \;,
\end{equation}
so $Y_i \sim \operatorname{Ber}(p)$. The behavior of $H_n$ at $\tilde t$ is uniquely determined by the $Y_i$, since
\begin{equation}
H_n(\tilde t)=- \ln \left( \tfrac 1n \sum_{i=1}^n Y_i\right) \; .
\end{equation}

\section{Establishing the Rate Functions}
\label{sec:Establishing}
\begin{figure}[!htb]
  \includegraphics[width=\textwidth]{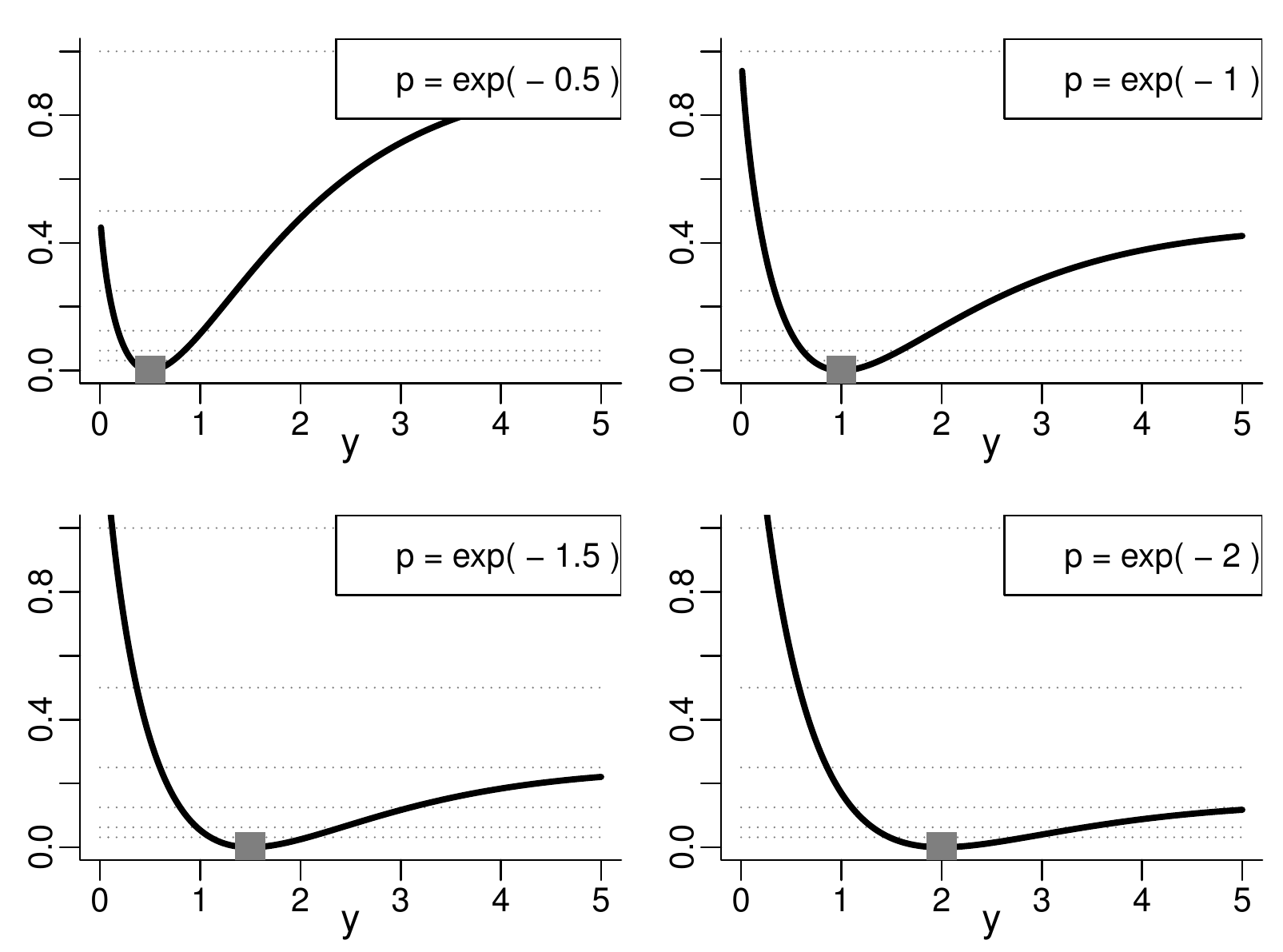}
\caption[b]{The rate function $I_p(y)$ (thick black line) from equation \eqref{eq:ldp2} for tail probabilities $p=\exp(-0.5)$ (top left), $p=\exp(-1)$ (top right), $p=\exp(-1.5)$ (bottom left), and $p=\exp(-2)$ (bottom right). The grey square is located at $- \ln (p)$, which is the minimum of the rate functions and the exact value of the cumulative hazard rate at the given tail probability. The dotted lines indicate the values $2^i$ for $i =0,-1,\dots,-5$.\\
	As the tail probabilities decrease, the rate function to the right of the grey square becomes increasingly flat, which indicates a low rate of convergence to the cumulative hazard rate. Although the same effect can be observed on the left side of the grey square, it is a lot weaker. The figure was generated using R\citep{R}.}
\label{fig:1}
\end{figure}

First, a  large deviation principle (LDP) for the $Y_i$ is derived. By either Cram\'{e}r's theorem \citep{cramer1938}\citep[p. 508]{klenke2008} or by Sanov's theorem \citep{sanov1958}\citep[p. 518]{klenke2008}, the series of probability measures
\begin{equation}
(P_n)_{n \in \mathbb N}=  \left( \mathcal L \left( \tfrac 1n \sum_{i=1}^n Y_i\right) \right)_{n \in \mathbb N}
\end{equation}
satisfies a large deviation principle with rate $n$ and rate function
\begin{equation}
\tilde I_p(x)=x \ln \left( \frac{x}{p}\right) + (1-x) \ln \left( \frac{1-x}{1-p}\right) \label{eq:ldp1}
\end{equation}
for $x \in [0,1]$. This is the relative entropy of two Bernoulli distributions, one with success probability $p$ and one with success probability $x$\citep[p. 515]{klenke2008}.\\
Applying the contraction principle \citep[p. 518]{klenke2008} to this LDP and the function $f(x)= - \ln(x)$
shows that the series $\left(\mathcal L \left(-\ln \left( \tfrac 1n \sum_{k=1}^n Y_k\right)\right)\right)_{n \in \mathbb N}$ satisfies a LDP with rate $n$ and rate function
\begin{align}
I_p(y)&:= \tilde I_p(f^{-1}(y))\\
&=\exp(-y) \ln \left( \frac{\exp(-y)}{p}\right) + (1-\exp(-y)) \ln \left( \frac{1-\exp(-y)}{1-p}\right) \label{eq:ldp2}
\end{align}
for $y \in [0, + \infty]$. This is the rate function of $H_n(\tilde t)$, and is displayed in fig. \ref{fig:1} for four different values of $p$.
Next, substitute $y=z-\ln(p)$ and define the centered rate function as 
\begin{align}J_p(z)&:=I_p(z-\ln(p)) \\
&=-z p \exp(-z) + (1-p\exp(-z)) \ln \left( \frac{1-p\exp(-z)}{1-p}\right)\label{eq:ldp3}
\end{align}
for $z \in [\ln(p),+\infty]$ and $p \in (0,1)$. It is defined since the main interest of this examination is to compare the behavior of the rate functions $I_p$ for different $p$ close to the pointwise limit value $-\ln(p)$ of the empirical cumulative hazard rate. This value is shifted to the origin in the centered rate function and therefore allows to compare the behavior of the empirical cumulative hazard rate at fixed distances from the respective limit values for different $p$.

\section{Properties of the Rate Functions}
\label{sec:Properties}
\subsection{Monotonicity in \textit{p}}
\label{sec:Monotonicity}
In this section, it is shown that $J_p$, taken as a function of $p$, is strictly increasing. This implies that the speed of convergence is decreasing as $\tilde t$ increases and correspondingly the tail probabilities decrease. Without loss of generality, it is assumed that $X$ takes on every tail probability, so $J_p$ is well defined for all $p \in (0,1)$.\\
First, the function $J_p$ can not be expected to be strictly increasing in $p$ for $z=0$, since by definition $J_p(0)=0$ for all $p\in (0,1)$. Therefore the case $z=0$ will be excluded.\\
The first derivative of $J_p$ with respect to $p$ is given by
\begin{equation}
\frac{ d}{ dp} J_p(z)=-z \exp(-z)  + \frac{1-\exp(-z)}{1-p} -\exp(-z) \ln \left(  \frac{1-p\exp(-z)}{1-p}\right)
\end{equation}
and the second derivative of $J_p$ with respect to $p$ by
\begin{equation}
\frac{  d^2}{  dp^2} J_p(z)=  \frac{\exp(-z)(\exp(z)-1)^2}{(1-p)^2(\exp(z)-p)}\; . \label{eq:mon}
\end{equation}
Since the inequality $\exp(z)-p >0$ holds for all $p \in (0,1)$ and all $z$ in the domain of $J_p$, termwise analysis of eq. \eqref{eq:mon} shows that the second derivative is positive for all feasible $p$, and zero only when $z=0$, which was excluded above. This makes the second derivative strictly positive, so $J_p$ is strictly convex in $p$.\\
If $J_p$ as defined in eq. \eqref{eq:ldp3} is taken as a function of $p$ for $p<1$, it is well-defined and its first derivative, evaluated at $p=0$, yields
\begin{equation}
\frac{d}{dp} J_0(z) = \exp(-z)(-z-1)+1 >0 \text{ for } z \neq 0 \; .
\end{equation}
So $J_p$ is strictly convex in $p$ and its gradient at $p=0$ is strictly positive, therefore $J_p$ is strictly increasing in $p$ for $p \in (0,1)$. This shows that the rate of convergence $J_p(z)$ is decreasing for all $z$ as the tail probabilities $p$ decrease.
\subsection{Asymmetry in \textit{z}}
\label{sec:Asymmetry}
First, it is shown that $J_p$ is not axis symmetric with respect to the ordinate axis. As an aid, the identity
\begin{equation}
(1-x)\ln(1-x)=-x+\sum_{k=2}^\infty \frac{x^k}{k^2-k} \;\text{ for } |x|<1\label{eq:aid}
\end{equation}
is used. Splitting the fraction in the logarithm in $J_p$ from equation \eqref{eq:ldp3} and then applying the power series from eq. \eqref{eq:aid} with $x=p\exp(-z)$ yields
\begin{equation}
J_p(z)=p\exp(-z) \left[-z +\ln \left( \tfrac{1-p}{e} \right)+ \sum_{k=2}^\infty \frac{(p\exp(-z))^{k-1}}{k^2-k}\right] - \ln(1-p) \label{eq:asy}
\end{equation}
for $p \exp(-z)<1$. Therefore $J_p$ is asymmetric by the asymmetry of $\exp(-z)$.\\
Next, the symmetry defect of $J_p$, given by $|J_p(-z)-J_p(z)|$, is examined. For this, the identities
\begin{equation}
2\sinh(z)=\exp(z)- \exp(-z) \text{ and } 2 \cosh(z)= \exp(z)+ \exp(-z)
\end{equation}
are used. With the representation of $J_p$ in eq. \eqref{eq:asy}, directly subtracting $J_p(z)$ and $J_p(-z)$ yields
\begin{align}
\big \vert J_p(-z)-J_p(z) \big \vert &= \Big \vert  \;2zp \cosh(z)+ 2 p\ln(\tfrac{1-p}{e})\sinh(z) \\
& \quad \quad + \sum_{k=2}^\infty \frac{2p^k\sinh(kz)}{k^2-k} \; \Big \vert\; .
\end{align}
Using only the second order term of the sum gives the approximation for the symmetry defect
\begin{equation}
\big \vert J_p(-z)-J_p(z) \big \vert \approx 2p \left(z\cosh(z) + \sinh(z)\ln(\tfrac{1-p}{e})\right) + p^2 \sinh(2z)
\end{equation}
for $p \in (0,1)$ and $z \in (0, -\ln(p))$.
These two statements show that the rate of convergence is not symmetric around the true cumulative hazard rate $H$, and will always be slower from above the cumulative hazard rate than from below.
\section{Example}
\label{sec:Example}
\begin{figure}[!htb]
  \includegraphics[width=\textwidth]{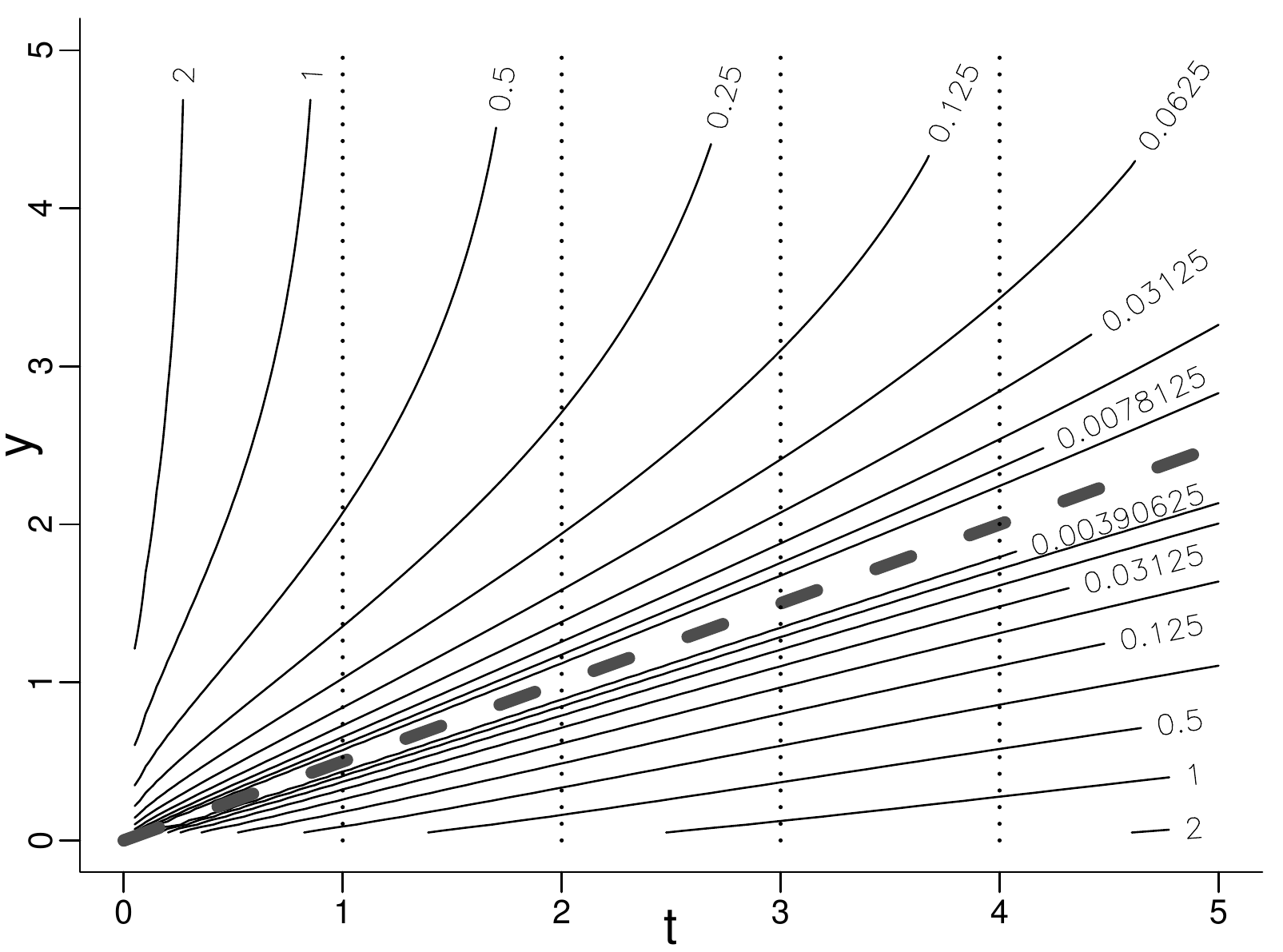}
\caption[b]{A contour plot showing lines of equal rate of convergence (thin black lines) to a cumulative hazard rate (thick grey dashed line). The black contour lines are determined by $I(y,t)$ from eq. \eqref{eq:twopara} for the special case of an exponential distribution with mean $\tfrac 1 \lambda=2$ and are lines of equal rate of convergence for the rates $2^i$, where $i=2,1,0,\dots,-8$. The grey dashed line is the cumulative hazard rate, which is a linear function through the origin with gradient $\tfrac 12$. The four rate functions displayed in fig. \ref{fig:1} are sections through this contour plot along the dotted black lines: The grey squares in fig. \ref{fig:1} are the intersection of the dotted black lines and the dashed grey cumulative hazard rate, the thin dashed lines in fig. \ref{fig:1} represent some of the contour lines from the figure shown here.\\
	As $t$ increases and the tail probabilities decrease accordingly, the lines on the upper left half of the picture diverge faster from the cumulative hazard rate than on the lower right half of the picture, showing that the rate of convergence is lower on the upper left side.\\
	Note that since the rate functions are only determined by the tail probabilities, and all tail probabilities occur in an exponential distribution, the corresponding contour plot for any random variable can be derived from this plot by a transformation of the $t$-axis. The figure was generated using R\citep{R}.}
\label{fig:2}
\end{figure}

Let $X$ be a positive random variable with distribution function $F$ and survival function $S$. Taking the tail probabilities $p$ as functions of $t$, meaning $p=p(t)=S(t)$, the rate function $I_p$ in equation \eqref{eq:ldp2} is a function in the variables $t$ and $y$:
\begin{align}
I(y,t)&:=I_{p(t)}(y)\\
&=\exp(-y) \ln \left( \frac{\exp(-y)}{S(t)}\right) + (1-\exp(-y)) \ln \left( \frac{1-\exp(-y)}{F(t)}\right). \label{eq:twopara}
\end{align}
For the case where $X$ is exponentially distributed with mean $\tfrac 1 \lambda =2$, the contour plot of $I(y,t)$ is displayed in fig. \ref{fig:2}, alongside the cumulative hazard rate. The asymmetry of the rate functions is clearly visible.

\section{Discussion}
\label{sec:Discussion}
The main result of this paper is that the estimate of the cumulative hazard rate (short CH) has a tendency to overestimate the true CH, which is strongest in the cases of (1) small sample sizes or (2) small tail probabilities. This is a direct conclusion from the asymmetry of the rate function governing the convergence of the estimated cumulative hazard rate in combination with the fundamental relation of the theory of large deviations given in equation \eqref{eq:ldpfundamental}. The degree of this effect is determined by the tail probabilities and therefore unique for every distribution.\\
It is notable that the asymmetry observed is not generated by the distribution of the random variables, but is rooted in the asymmetry of the rate function of the Bernoulli distribution. It is amplified by the logarithm in the definition of the cumulative hazard rate, and unavoidable in the sense that it naturally arises from observing whether a random variable takes on a value over or under a given threshold.\\
What remains unclear is how these results translate into the standard framework of survival analysis, i.e. censoring and using point processes instead of i.i.d. random variables. Especially censoring might have a strong influence on the case with low tail probabilities, since early censoring can render long survival times (that are commonly associated with low tail probabilities) irrelevant. Extending the results presented in this papers to this more general setting is relevant for applications, but requires further work.

\end{document}